\documentclass[fleqn,12pt]{article}
\usepackage[english]{babel}

\usepackage{amssymb,amsmath}
\newtheorem{Theorem}{Theorem}[section]
\newtheorem{Proposition}[Theorem]{Proposition}

\newtheorem{Lemma}[Theorem]{Lemma}

\newtheorem{Ex}{Example}[section]
\newtheorem{proofhead}{Proof}

\newenvironment{proof}
    {
        \par
        \begin{proofhead}
        \normalfont
    }
    {   \qed
        \end{proofhead}
        \par
    }

\newcommand{\bb}{\mathbb}

\def \mean{{\bb E}}

\def \var{{\bb V}{\rm ar}}

\newcommand{\qed}{\hfill$\Box$}

\newcommand{\refs}[1]{(\ref{#1})}

\newcommand{\halmos}{\hfill $\Box$}

\newcounter{mylistcnt}
\renewcommand{\themylistcnt}{{\rm({\roman{mylistcnt}})}}

\newcounter{zad}

\begin{document}

\title{Estimates for moments of supremum of reflected fractional Brownian
motion}
\author{
Krzysztof D\c{e}bicki\thanks{This work was supported by MNiSW
Research Grant N N201 394137 (2009-2011)}, Agata Tomanek
\vspace*{.08in} \\
Mathematical Institute, University of Wroc\l aw \\
pl. Grunwaldzki 2/4, 50-384 Wroc\l aw, Poland }

\date{}
\maketitle

\begin{abstract}
Let $B_H(\cdot)$ be a fractional Brownian motion with Hurst parameter $H\in(0,1]$.
Motivated by applications
to maximal inequalities for fractional Brownian motion,
in this note we derive bounds for
\[
K_T(H,\gamma):=\mean\left[\sup_{t\in[0,T]}|B_H(t)|\right]^\gamma,
\]
with $\gamma, T>0$.


\end{abstract}

\section{Introduction}

In this note we study properties of
\begin{eqnarray}
K_T(H,\gamma):=\mean\left[\sup_{t\in[0,T]}|B_H(t)|\right]^\gamma
\end{eqnarray}
for $\gamma, T>0$, where
$B_H(\cdot)$ is a fractional Brownian motion with Hurst parameter $H\in(0,1]$,
i.e., a centered Gaussian process with stationary increments and variance function
$\var (B_H(t))=t^{2H}$.

Constants $K_T(H,\gamma)$ appear in the context of analyzing maximal
inequalities for fractional Brownian motion; see e.g.
\cite{NoV99}.
The aim of this paper is to give bounds for
$K_T(H,\gamma)$.
To our best knowledge the exact values of $K_T(H,\gamma)$ are not known.

Let
\[
\beta(s)=\sum_{k=0}^\infty (-1)^k(2k+1)^{-s},
\]
$s>0$, be the {\it Dirichlet beta function}. By
$E_n(x)$, $n=0,1,...$, we denote the Euler polynomials;
$E_n:=2^n E_n(1/2)$, $n=0,1,...$, stand for Euler numbers; see \cite{AbS}, page 804.
Additionally, let $\Psi(t)=\mathbb{P} (\mathcal{N}>t)$ where
$\mathcal{N}\sim \mathcal{N}( 0,1)$.

The use of technique based on comparison of Gaussian processes yields the following  theorem.

\begin{Theorem}\label{th.2}
Let $\gamma>0$.
\\
(i) If $H<1/2$, then
$$K_T(H,\gamma)\ge T^{\gamma H} \frac{1}{\sqrt{\pi}} 2^{\frac{\gamma}{2}} \Gamma\left( \frac{\gamma+1}{2} \right).$$
\\
(ii) If $H\ge 1/2$, then
$$T^{\gamma H} \frac{1}{\sqrt{\pi}} 2^\frac{\gamma}{2}\Gamma\left( \frac{\gamma+1}{2} \right)\le K_T(H,\gamma)\le T^{\gamma H} \frac{1}{\sqrt{\pi}} 2^{\frac{\gamma}{2}+1} \Gamma\left( \frac{\gamma+1}{2} \right).$$
\end{Theorem}

In the following proposition we calculate exact value of $K_T(1/2,\gamma)$ and $K_T(1,\gamma)$.

\begin{Proposition}\label{th.1}
Let $\gamma>0$. Then
\\
(i)
$K_T(1/2,\gamma)=\frac{1}{\sqrt{\pi}}2^{1+\frac{\gamma}{2}}\Gamma\left(
\frac{\gamma+1}{2} \right)\beta(\gamma)T^\frac{\gamma}{2}$;
\\
(ii)
$K_T(1,\gamma)=\frac{1}{\sqrt{\pi}}2^{\frac{\gamma}{2}}\Gamma\left(
\frac{\gamma+1}{2} \right)T^\gamma$.
\end{Proposition}
The detailed proofs of Proposition \ref{th.1} and Theorem \ref{th.2} are deferred to Section \ref{s.proofs}.

As an immediate consequence of Proposition \ref{th.1}, in view of \cite{AbS}, page 805, we have
\begin{eqnarray}
K_1(1/2,2n+1)=\sqrt{\frac{\pi}{2}} \left( \frac{\pi^2}{2}\right)^n
\frac{n!}{(2n)!}|E_{2n}|\label{eq.odd}
\end{eqnarray}
\[
K_1(1/2,2n)=\frac{(-1)^n}{(n-1)!}\left( \frac{\pi^2}{2}\right)^n\int_0^1 E_{2n-1}(x)sec(\pi x)dx,
\]
for $n=1,2,...$.
In particular $K_1(1/2,1)=\sqrt{\pi/2}$ and
$K_{1/2}(1/2,2)$ is the {\it Catalan's constant}.

The comparison of the upper bound for $K_1(1/2,2n+1)$ given in Theorem \ref{th.2}
with \refs{eq.odd} enables us to recover
the known inequality for Euler numbers
\[
|E_{2n}|\le\frac{4^{n+1} (2n)!}{\pi^{2n+1}};
\]
see, e.g., \cite{AbS}, page 805.
We refer to \cite{DaG02} for other results that relate moments of functionals of Brownian motion with number theory.

\section{Proofs}\label{s.proofs}
In this section we present complete proofs of Proposition \ref{th.1} and Theorem \ref{th.2}.

We frequently use the fact that
the property of self-similarity of fractional Brownian motion
enables us to write
\begin{eqnarray}
K_T(H,\gamma)=K_1(H,\gamma) T^{\gamma H}.\label{reduction}
\end{eqnarray}

We start with an auxiliary result which is also of independent interest.
\begin{Lemma}\label{l.gen}
Let $\{X(t):t\ge0\}$ be a centered Gaussian process with stationary increments and continuous and strictly increasing variance function $\sigma^2_X(\cdot)$,
$X(0)=0$ a.s.
\\
(i) If $\sigma^2_X(\cdot)$ is concave, then
\[
\mean\left[\sup_{t\in[0,T]}X(t)\right]^\gamma\ge \left(\sigma^2_X(T)\right)^\frac{\gamma}{2} \frac{1}{\sqrt{\pi}} 2^\frac{\gamma}{2} \Gamma\left( \frac{\gamma+1}{2} \right) ;
\]
(ii) If $\sigma^2_X(\cdot)$ is convex, then
\[
\mean\left[\sup_{t\in[0,T]}X(t)\right]^\gamma\le \left(\sigma^2_X(T)\right)^\frac{\gamma}{2} \frac{1}{\sqrt{\pi}} 2^\frac{\gamma}{2} \Gamma\left( \frac{\gamma+1}{2} \right)  .
\]
\end{Lemma}
\proof
Since the proof of $(ii)$ is analogous to the proof of $(i)$, we focus on the argument that justifies
$(i)$.
Assume that $\sigma^2_X(\cdot)$ is concave.
Observe that for $Y(t):=B_{\frac{1}{2}}\left(\sigma_X^2(t)\right)$ we have
$$\mathbb{V}ar\left(Y(t)\right)=\mathbb{V}ar\left(B_{\frac{1}{2}}\left(\sigma_X^2(t)\right)\right)=\sigma_X^2(t)=\mathbb{V}ar(X(t))$$
for all $t\in[0,T]$ and, due to concavity of $\sigma_X^2(\cdot)$,
\begin{eqnarray}
\nonumber \mathbb{V}ar\left(Y(t)-Y(s)\right) & = & \mathbb{V}ar\left(B_{\frac{1}{2}}\left(\sigma_X^2(t)\right)-B_{\frac{1}{2}}\left(\sigma_X^2(s)\right)\right)\\
\nonumber & = & \sigma_X^2(t)-\sigma_X^2(s)\\
&\le& \sigma_X^2(t-s) \nonumber \\
\nonumber & = & \mathbb{V}ar(X(t)-X(s))
\end{eqnarray}
for all $t>s$ and $s,t\in[0,T]$.
Thus, using Slepian inequality (see, e.g., Theorem 2.1 in Adler \cite{Adl90}),
$$ \mathbb{P}\left(\sup_{t\in [0,T]}X(t)>x \right)\ge \mathbb{P}\left(\sup_{t\in [0,T]}B_{\frac{1}{2}}(\sigma_X^2(t))>x \right)$$
for all $x\ge 0$.
Since $\mathbb{P}\left(\sup_{t\in [0,T]}B_{\frac{1}{2}}\left(\sigma_X^2(t)\right)>x \right) =  \mathbb{P}\left(\sup_{t\in\left[0,\sigma_X^2(T)\right]}B_{\frac{1}{2}}(t)>x \right)$,
then we get
$$\mathbb{E}\left[\sup_{t\in[0,T]}X(t) \right]^\gamma\ge \  \mathbb{E}\left[ \sup_{t\in\left[0, \sigma_X^2(T)\right]}B_{\frac{1}{2}}(t) \right]^\gamma. $$
Due to self-similarity of Brownian
motion we have
\begin{eqnarray}
\nonumber \mathbb{E}\left[ \sup_{t\in\left[0,\sigma_X^2(T)\right]}B_{\frac{1}{2}}(t) \right]^\gamma & = & \mathbb{E}\left[ \sup_{t\in[0,1]}B_{\frac{1}{2}}\left( \sigma_X^2(T)\ t\right) \right]^\gamma\\
& = & \left( \sigma_X^2(T)\right)^\frac{\gamma}{2} \mathbb{E}\left[ \sup_{t\in[0,1]}B_{\frac{1}{2}}(t) \right]^\gamma.
\label{in.l}
\end{eqnarray}
Finally, using that
$\mathbb{P}\left(\sup_{t\in[0,1]}B_{\frac{1}{2}}(t)>t\right)=2
\mathbb{P}\left(\mathcal{N}>t\right)$,
we get
\begin{eqnarray}
\mathbb{E}\left[ \sup_{t\in[0,1]}B_{\frac{1}{2}}(t) \right]^\gamma=
\int_0^\infty \gamma x^{\gamma-1}2\Psi\left( x\right)dx =  \frac{1}{\sqrt{\pi}} 2^\frac{\gamma}{2}\Gamma\left( \frac{\gamma+1}{2} \right).\label{for.last}
\end{eqnarray}
Combination of \refs{in.l} with
\refs{for.last}
completes the proof of $(i)$.
\endproof

\subsection{Proof of Theorem \ref{th.2}}
Following \refs{reduction} we consider only the case of $T=1$.
Note that
$$\mathbb{P}\left(\sup_{t\in[0,1]}B_H(t)>x \right)  \le \mathbb{P}\left(\sup_{t\in[0,1]}|B_H(t)|>x \right) \le 2\  \mathbb{P}\left(\sup_{t\in[0,1]}B_H(t)>x \right). $$
The combination of the above with Lemma \ref{l.gen} completes the proof of $(i)$ and the upper bound in $(ii)$.\\
To prove the lower bound in $(ii)$ we use that
\begin{eqnarray*}
 \mathbb{E}\left[ \sup_{t\in\left[0, 1\right]}|B_H(t)| \right]^\gamma
 &\ge&
 \mathbb{E} |B_H(1)|^\gamma=
\mathbb{E} |\mathcal{N}|^\gamma=
\frac{1}{\sqrt{\pi}} 2^\frac{\gamma}{2}\Gamma\left( \frac{\gamma+1}{2} \right).
\end{eqnarray*}
This completes the proof.

\halmos

\subsection{Proof of Proposition \ref{th.1}}

{\bf Ad $(i).$}
Assume that $H=\frac{1}{2}$.
Following \refs{reduction} it suffices to analyze
$K_1\left(1/2,\gamma\right)$.
Using formula 1.1.4 in \cite{BorSal}, we
get
\begin{eqnarray}
K_1\left(1/2,\gamma\right)
&=& \mathbb{E}\left[ \sup_{t\in[0,1]} \left|B_{\frac{1}{2}}(t)\right| \right]^\gamma \nonumber\\
&=&\gamma \int_0^{\infty} x^{\gamma-1} \mathbb{P}\left( \sup_{t\in [0,1]} \left|B_{\frac{1}{2}}(t)\right|>x \right)dx\nonumber\\
&=&
\gamma \int_0^{\infty} x^{\gamma-1}\sum_{k=-\infty}^{\infty}\left\{ (-1)^k sign((2k+1)x)\frac{2}{\sqrt{\pi}}\int_{\frac{x|2k+1|}{\sqrt{2}}}^{\infty} e^{-s^2}ds\right\} dx.\nonumber\\
\label{sum.1}
\end{eqnarray}
Let
\begin{eqnarray*}
f_n(x)&:=&x^{\gamma-1}\sum_{k=-n-1}^{n}\left\{ (-1)^k sign((2k+1)x)\frac{2}{\sqrt{\pi}}\int_{\frac{x|2k+1|}{\sqrt{2}}}^{\infty} e^{-s^2}ds\right\}\\
&=&
2
x^{\gamma-1}\sum_{k=0}^{n}\left\{ (-1)^k \frac{2}{\sqrt{\pi}}\int_{\frac{x|2k+1|}{\sqrt{2}}}^{\infty} e^{-s^2}ds\right\}
\end{eqnarray*}
and observe that
for each $n\ge 0$ and $x\in(0,1]$
\begin{eqnarray}
0\le f_n(x)\le f_0(x)\le 2 x^{\gamma-1}.\label{leb1}
\end{eqnarray}
Additionally, for each $n\ge 0$ and $x>1$, we have
\begin{eqnarray}
\nonumber |f_n(x)| & \leq &  x^{\gamma-1}4\sum_{k=0}^{n}\frac{1}{\sqrt{2\pi}}\int_{x(2k+1)}^\infty e^{-\frac{s^2}{2}}ds=4x^{\gamma-1}\sum_{k=0}^{n}\Psi(x(2k+1))\\
\nonumber & \leq & \frac{4}{\sqrt{2\pi}}x^{\gamma-1}\sum_{k=0}^{n}\frac{1}{x(2k+1)}e^{-\frac{x^2(2k+1)^2}{2}}
\label{upper.bound.g}\\
& \leq & \frac{4}{\sqrt{2\pi}}x^{\gamma-2}e^{-\frac{x^2}{2}}\sum_{k=0}^\infty \left( e^{-x^2} \right)^k \\
 & = & \frac{4}{\sqrt{2\pi}}x^{\gamma-2}e^{-\frac{x^2}{2}} \frac{1}{1-e^{-x^2}},\label{leb2}
\end{eqnarray}
where \refs{upper.bound.g} follows from the fact that
$\Psi(t)\le\frac{1}{\sqrt{2\pi}t}\exp(-\frac{t^2}{2})$ for each $t\ge 0$.
The combination of \refs{leb1} with \refs{leb2}
implies that $|f_n(\cdot)|$ is bounded by an integrable function, and
hence, by Lebesgue's dominated convergence theorem, we can rewrite \refs{sum.1}
in the following form
\begin{eqnarray}
\nonumber K_1\left(1/2,\gamma\right)
\nonumber & = & \gamma \frac{4}{\sqrt{\pi}} \sum_{k=0}^\infty \left\{ (-1)^k \underbrace{\int_0^\infty \left\{ x^{\gamma-1} \int_{\frac{x(2k+1)}{\sqrt{2}}}^\infty e^{-s^2} ds \right\}dx}_{I_k} \right\}.
\end{eqnarray}
The change of the order of integration in $I_k$ leads to
\begin{eqnarray}
\nonumber I_k
\nonumber & = & \frac{1}{\gamma} \left( \frac{\sqrt{2}}{2k+1} \right)^\gamma \int_0^\infty s^\gamma e^{-s^2}ds=\frac{1}{2\gamma} \left( \frac{\sqrt{2}}{2k+1} \right)^\gamma \int_0^\infty e^{-t} t^{\frac{\gamma-1}{2}}dt \\
\nonumber & = & \frac{1}{2\gamma} \left( \frac{\sqrt{2}}{2k+1} \right)^\gamma \Gamma\left( \frac{\gamma+1}{2} \right),
\end{eqnarray}
which implies that
\begin{eqnarray}
K_1\left(1/2,\gamma\right)
&=&
\frac{1}{\sqrt{\pi}}2^{\frac{\gamma}{2}+1} \Gamma\left( \frac{\gamma+1}{2} \right)
\sum_{k=0}^\infty (-1)^k \left( \frac{1}{2k+1} \right)^\gamma
=
\frac{1}{\sqrt{\pi}}2^{\frac{\gamma}{2}+1}\Gamma\left( \frac{\gamma+1}{2} \right) \beta(\gamma)\nonumber.
\end{eqnarray}
This completes the proof of $(i)$.
\\
\\
{\bf Ad $(ii)$.}
Let $H=1$. Since $B_{1}(t)=_{\rm d}\mathcal{N}t$, where $\mathcal{N}\sim \mathcal{N}(0,1)$, then
$\mathbb{P}\left( \sup_{t\in [0,1]}|B_1(t)| >x\right)= 2\Psi(x)$.
Standard integration completes the proof.
\halmos


\begin{thebibliography}{99}\small


\bibitem{AbS} Abramowitz, M., Stegun, I.A. {\it Handbook of Mathematical Functions
with
Formulas, Graphs, and Mathematical Tables}, National Bureau of Standards
Applied Mathematics Series - 55, Washington, 1972.
 \bibitem{Adl90}
 Adler, R.J. \emph{An introduction to continuity, extrema, and related topics for general Gaussian processes}
 Inst. Math. Statist. Lecture Notes -Monograph Series, vol. 12, Inst. Math. Statist., Hayward, CA, 1990.
 \bibitem{BorSal}
 Borodin, A.N., Salminen, P. \emph{Handbook of Brownian Motion - Facts and Formulae} Birkh\"{a}user Verlag, Basel, Boston, Berlin, 1996.
\bibitem {DaG02}
DasGupta, A. (2002) Mellin transform and densities of suprema of brownian processes: with applications.
Report 02-08, Department of Statistics, Purdue University, USA.

\bibitem{NoV99}
Novikov, A., Valkeila, E. (1999) On some maximal inequalities for fractional Brownian motions.
{\it Statistics and Probability Letters} {\bf 44}, 47--54.


%
%
%
%
\end{thebibliography}
\end{document}